\documentclass[hein,final,leqno]{shein}
\usepackage{latexsym}
\usepackage{graphics}
\usepackage{epsfig}
\usepackage{epsf}
\usepackage{eucal}
\usepackage{pstricks}
\usepackage{amsfonts}
\usepackage{amstext}
\usepackage{amssymb}
\usepackage{amsmath}
\usepackage{amsbsy}
\usepackage{amsxtra}
\usepackage{amscd}
\shnewtheorem*{udefinition}{Definition}{\bfseries\itshape}{} 
\shnewtheorem*{uremark}{Remark}{\bfseries\itshape}{}   
\shnewtheorem*{uremarks}{Remarks}{\bfseries\itshape}{} 
\shnewtheorem*{utheorem}{Theorem.}{\bfseries}{} 
\shnewtheorem*{uproposition}{Proposition}{\bfseries}{} 
\shnewtheorem*{ucorollary}{Corollary.}{\bfseries}{} 
\renewcommand{\abstractname}{Abstract.}
\title {Gauge techniques in time and frequency domain TLM}
\author{Steffen Hein\\ \\\small{Spinner GmbH. M\"unchen\\Aiblinger Str. 30,
DE-83620 Feldkirchen-Westerham, Germany }}
\renewcommand{\abstractname}{Abstract.}
\begin{document}
\institute{Steffen Hein 
\hfill DE-83043 Bad Aibling, Germany\hfill}
\titlerunning{Gauge techniques in time and frequency domain TLM}
\authorrunning{Steffen Hein}
\vspace{-3cm}
\maketitle 
%
\begin{abstract}
Typical features of the Transmission Line Matrix (TLM) algorithm in connection
with stub loading techniques and prone to be hidden by common frequency domain
formulations are elucidated within the propagator approach. In particular, the
latter reflects properly the perturbative character of the TLM scheme and
its relation to gauge field models. Internal `gauge' degrees of freedom are
made explicit in the frequency domain by introducing the complex nodal S-matrix
as a function of operators that act on external or internal fields or virtually
couple the two. As a main benefit, many techniques and results gained in the
time domain thus generalize straight away. The recently developed deflection
method for algorithm synthesis, which is extended in this paper, or the
non-orthogonal node approximating Maxwell's equations, for instance, become
so at once available in the frequency domain. In view of applications in
computational plasma physics, the TLM model of a relativistic charged particle
current coupled to the Maxwell field is traced out as a prototype.
{\vspace{2pt}\hfill
\textbf{MSC-classes}:\;\textnormal{65C20,\,65M06,\,76D05}}
\medskip
\vspace{-10pt}
\end{abstract}
%
\markboth{{\normalsize \textsc{Steffen Hein}}}
{{\normalsize \textsc{Time and frequency domain TLM}}}
\normalsize
\vspace{-19pt}
\section {Introduction}
When P. B. Johns and coworkers introduced the transmission line matrix (TLM)
numerical method in the early seventies \cite{John5}  they were certainly not
aware of handling gauge field concepts. In fact, the perturbative aggregation
of internal gauge degrees of freedom plays a salient role in TLM in the form
of  `stub loading'. This well known technique, ad hoc applied  in a great many
variants to model continuous mesh deformations or material parameters, e.g.,
expands to a powerful framework that allows for a systematic integration of
manifold interactions. Stubs in a generalized sense may for instance add mass
to a field much alike a Higgs field \cite{Hein12}, or induce superconducting
or non-reciprocal behaviour \cite{Hein10},\cite{Hein11}. Very much of the
strength and versatility of TLM rests on this feature that indeed relates the
method to some pioneering developments in modern particle physics. 

Varied formal characterizations of the TLM scheme have been given since the
early transmission line network interpretations. Some authors proceed
from Huygens' principle \cite{Hoef} (sometimes misunderstood) others use
concepts of state space control theory \cite{Wit7} or random walk
\cite{Cogan6}, \cite{Cogan7}. The present paper is not intended to add a
further characterization to these and alternate
interpretations which are usually built on the cinematics of a special kind
of propagating quantities (waves, particles, temperatures etc.). Instead, this
study focusses on the internal structure of the TLM algorithm irrespective of
its physical interpretation. It nevertheless remains remarkable that so
many heterogeneous perceptions exist and apparently match the algorithmic
scheme of the method. The perturbative structure of TLM, with its gauge
fixture not at least, certainly explains much of that multiplicity.

With sustained evolution of computer performance in coming years there will
be an increasing need of highly flexible numerical schemes, complementary to
powerful standard solvers. The versatile structure and the diversity of yet
known applications in wave propagation \cite{Chris1}, chemical
reaction-diffusion  \cite{Wit7},
heat transfer  \cite{Cogan7}, plasma physics \cite{Hein12}, and further
domains provide some evidence that the TLM
method can play a competent role on the court.

This requires, however, a refined  technical outfit being built on sound
mathematics. Theorems like the deflection lemma \cite{Hein11}, a generalized
version of which is given in this paper, should constitute a canonical basis for
algorithm synthesis. The usage and usefulness of those tools amenable to the
propagator approach  \cite{Hein8} have already been shown up within a series of
numerical applications  \cite{Hein9} -  \cite{Hein12}. On the issue of this 
paper the TLM model of a relativistic charged particle current coupled to the 
Maxwell field provides supplementary illustration.

\vspace{-8pt}
\section {The algorithm and its interpretation}
This section is widely  heuristic and starts slightly outside our subject.
Conceptually, the TLM algorithm will be separated from its interpretation in
terms of any modelled fields. This obviously implies outlining a canonical
passage between two severed parts: physical fields and computed quantities.
The TLM typical manner of generating these quantities, i.e. the computational
part, is then headed for necessarily on a  level of abstraction beyond
the commonly considered Maxwell field.

It is in the obvious nature of any splittings (conceptual just as physical)
that junctions appear double faced, and so the joints between the TLM
algorithm and its interpretation, viz. ports and nodes in cells, take a
twofold role. So it has to be clarified what is meant in either context.

At the interface to  physical fields in space a {\it  cell} stands for a
convex polyhedron in configuration space ${\mathbb R}^D$ and a {\it port} for a
distribution, i.e. a continuous linear functional acting on a field, with its
support on a closed subset of the cell boundary (a point, line segment, or face
e.g.). The distributional support is also referred to as the {\it location} 
of that port. If no confusion can arise it will be simply identified with
the port.

Again, in the same context, the {\it node} comprises the family of all port
distributions of a cell, yet shifted by  translations into its interior (where
conceptually one effective internal state of the cell at a time interacts with
boundary states - this touches the algorithmic role of ports and nodes treated
below). Note that in many cases distributed integrals instead of point
supported functionals, viz. the shifted {\it nodal images} of the ports, probe
the fields in the cells.

In a condensed slightly formalized writing: Every port $p$ of a cell together
with the node gives rise to a pair of linear forms $(p,Z)$, $(p,Z)^{\sim}$
which probe a real or complex vector field $Z = (Z_{1} ,\ldots, Z_{m})$; $m \in
\mathbb{N}_{+}$, at the cell boundary and within the cell. The two forms are
naturally interrelated by pull-back
\vspace{-2pt}
\begin{equation}\label{2.1} \left( \; p, Z \;
\right)^{\sim} \quad = \quad \left( \; p \circ s, Z \; \right) \quad = \quad
\left( \; p, Z \circ  s^{-1} \right)
\end{equation}
where  $s$  denotes the pertinent translational shift of the port onto its
nodal image, cf. fig 1.
Defining these forms consistently as distributed or point supported
functionals ({\it finite integrals}, e.g.) is a constituent part of any
interpretation of the TLM model.

\begin{figure}[h] 
\input{f1.tex}
\end{figure}

Quite different, and still more abstract,  are the roles of ports, nodes, and
cells in the computational context. Every cell stands here for a direct sum
(viz. a triple) of linear spaces, which respectively represent observable
{\it link} states in the cell or on its boundary, and unobservable {\it gauge}
(or {\it stub}) degrees of freedom. Observable vectors $
z_{\mu}^{p,n}$  are labelled by ports $p$ or their nodal images $n$, and are
linearly interrelated to physical fields by an  interpretation. In a more
precise diction, an {\it interpretation} may be defined as a family of linear
maps $I_\zeta$, indexed by the cells of a {\it mesh}, with domain  a distinguished 
class of vector fields in configuration space  which represent physical fields and
with their range in the set of observable cell states:
\begin{equation}\label{2.2} I_\zeta : Z_\mu \to \left( z_\mu ^{p,n} \right)
\; ; \;\; z_\mu^p = \left( \; p_\mu , Z_\mu \right), \quad
z_\mu^n=\left( \; p_\mu , Z_\mu \right)^\sim 
\end{equation} 
($\mu$ labels the ports and pertinent fields in cell $\zeta$ and superscripts
$p$, $n$ refer respectively to
ports on the cell boundary or to their nodal images in the cell - a port
being counted here as many times as it acts on a different field.) Note that
no fixed linear relations of such a kind interconnect any gauge vectors with
fields.

Ultimately, and by it most markedly characterizing the TLM scheme, the
ports together with their nodal images represent scattering channels for
states incident at the cell boundary and scattered in the node. As a
scattering channel every common port of two adjacent cells thus naturally
connects their nodes, i.e. constitutes a {\it link}.

The inherent scattering scheme of the algorithm postulates in the line of 
Lax-Phillips \cite{Lax2} a direct splitting of the link states into incoming
and outgoing components
\begin{equation}\label{2.3} z=z_{in} \oplus z_{out}
\end{equation}
which of course has to be defined in  harmony with the propagation properties
of the modelled fields. In reference \cite{Hein11} the details are carried through
for the Maxwell field, for example. The Poynting vector evaluated at the cell
faces yields there a canonical splitting and a quadratic form that is positive
or negative, respectively, on the incoming or outgoing subspaces.

In TLM diction the latter are usually  called the (spaces of) incident and 
reflected quantities and in the following $\pi _{in}$, $\pi_{out}$ denote the
pertinent projections.

The TLM scattering algorithm, which  is eo ipso perturbative,
solves finite difference equations in the form of model equations  between
states that in any interpretation should approximate physical fields (i.e. 
solutions of integro-differential equations) in any sense and order
(sect.4). This is the requirement of {\it consistence}, as a cardinal property
of the interpretation.

In contrast, the question of internal  convergence, which is essentially that
of computational {\it stability}, refers to intrinsic properties of the TLM
algorithm such as are subject of this paper. That question has thus to be
taken up in due course. (We are not unhappy that it is in general easier.)

\vspace{-8pt}
\section {The propagator relations}
Causality requires in the time domain the scattered  states of a cell to be
functions of previously incident states. With a cell dependent {\it reflection
operator} $\mathcal{R}$ the outgoing port quantities are thus 
\begin{equation}\label{3.1}
z_{out}^p\left(t + \tau \right) = \mathcal{R}\left(z_{in}^p \left(t- \mu \tau \right)
\right)_{\mu \in \mathbb{N}} 
\end{equation} 
where $\tau$ denotes the time step. (Suitable definitions release from introducing line impedances, such as are commonly used instead of the projections (\ref{2.3}) - in fact they are substitutes for projections, cf. the appendix.) $\mathcal{R}$ assembles the physical interactions  within the cell in a geometry dependent manner, and represents in this sense the TLM system {\it locally}, while the {\it global} topological structure, as given by the port linking
scheme and the boundary conditions, are brought in by the {\it connection map} $\mathcal{C}$. At every time step the latter transmits the reflected port quantities  to
adjacent cells or, at a system boundary, back into the same cell and
contingently adds the field excitations by acting linearly as 
\begin{equation}\label{3.2} \left(z_{in}^p \right)_{\zeta} =  \mathcal{C}
\left(\left(z_{out}^p \right) _\eta  ,   \left(z_{exc}^p \right)_\vartheta
\right) 
\end{equation} 
 where $\zeta$, $\eta$, $\vartheta$ label the cells and $z_{exc}^p$  denotes any excitation.

Short-cut, every pair ($\mathcal{C},\mathcal{R}_{\zeta}$)  consisting of a connection map $\mathcal{C}$
and a family $\mathcal{R}_{\zeta}$ of reflection operators labelled by the cells of the
same mesh will be referred to as a {\it TLM system}.

The interesting objects in a TLM system  are, of course, the causal operators
$\mathcal{R}_{\zeta}$  which like the pages of a book are merely bound together by the
`book binding' function $\mathcal{C}$. In addition to being defined on back-in-time
running sequences of incident states, $\mathcal{R}$  may itself depend on time and
may be non-linear. Before venturing on the TLM typical structure of 
$\mathcal{R}$ a glance is cast at some cinematical relations between port and node
quantities.

Operationally, the TLM approach, as a  genuine scattering scheme, splits off
a cinematical part from the dynamical evolution by handling some interactions
as a perturbation. In  classical Lax-Phillips scattering, for instance, two
`intertwingled' one parameter unitary groups describe respectively the free
and interacting dynamics \cite{Lax2}.  A related ansatz leads in perturbative quantum  field theory to the Dyson
expansion of S-parameters in the form of Feynman diagrams which again couple free
fields \cite{Mand4}. Technically, the TLM algorithm operates similarly with
freely propagating port and node quantities that before and past reflection
are forced into fixed phase relations by imposing in the time domain
\begin{equation}\label{3.3} {z_{in}^n\left(t \right) :=z_{in}^p \; ( \;t -
\frac{\tau}{2} \; ) \; \quad },\quad{z_{out}^n\left(t \right):=z_{out}^p \; ( \; t + \frac{\tau}{2} \; ) \; \quad } 
\end{equation}
Note well that these are technical arrangements to  enable the scattering
scheme, while only total quantities, 
\begin{equation}\label{3.4}
z^{p,n}=z_{in}^{p,n} \oplus
z_{out}^{p,n} \quad ,\qquad cf. \quad ( \ref {2.3})
\end{equation} 
represent observable states and in fact enter the model equations (cf. sect.4).

Scarce and innocent as they appear,  settings (\ref{3.1}, \ref{3.2}, \ref
{3.3}) together with (\ref {2.3}) contain yet very much of the substance of TLM.
On the basis of essentially these ingredients, which form the core
issues of the so-called propagator approach \cite{Hein8}, \cite{Hein12},
internal `gauge' degrees of freedom arise naturally on introducing additional
interactions by a perturbation, along with S-parameters that couple these  to
the original states. A simple way of giving a rigorous meaning to that is the
deflection lemma. Originally introduced in reference \cite{Hein11}, the lemma
is subsequently generalized to handle higher order equations.

\vspace{-8pt}
\section {The deflection lemma}
Finite difference equations that in the form of {\it model equations}
approximate  the dynamical equations of any physical interpretation
interrelate in the time domain sequences of total node and port quantities
\begin{equation}\label{4.1} \left[ z^{n,p} \right] = \left(z^{n,p}\left(t- \mu
\tau \right) \right)_{\mu =0,1,2,\ldots}  \end{equation} 
which represent observable states in the cell and on its boundary,  at present time $t$ and
in their history. A maximum $\mu$ that eventually enters the equations
obviously reflects the temporal order of the modelled dynamical problem, which
may be finite or infinite.  To {\it any} order the model equations can be
written 
\begin{equation}\label{4.2} \mathcal{F} [ z_{+}^n ]  [ z^p ] = 0 \end{equation} 
with a suitable causal operator $\mathcal{F}$ (where {\it causal} stands always for
being defined on back-in-time running sequences such as (\ref {4.1})).
subscript $+ (-)$ denotes a $\tau / 2$ time shift  
\begin{equation}\label{4.3}
[z_{+\left(- \right)}^{n,p} ]  := \left[z^{n,p} \right]_{t +
(- ) \frac{\tau}{2}} =\left(z^{n,p}\left( t + (- )\frac{\tau}{2} - \mu \tau \right)
\right)_{\mu \in \mathbb{N}} \end{equation} 
which in harmony with the phase relations
(\ref {3.3}) synchronizes port and node switching in equations (\ref {4.2}),
for every vector $z_{in}^p ( t )$  incident at a cell and constant
on each time interval $[ k \tau, ( k+1 ) \tau )$;  $ k\in \mathbb{Z}$. Note that $\mathcal{F}$ may itself be time dependent in like manner, i.e. as a time step function matching integer multiples of $\tau$.

In the literature equations (\ref {4.2}) are encountered  in manyfold guises.
For  instance, in reference \cite{Hein8} Maxwell's integral equations
discretized in a convex hexahedral cell of else arbitrary shape take the form
\begin{equation}\label{4.4} \psi_0 z^p \left( t \right) = \varphi_0  z^n
\; ( \; t + \frac{\tau}{2} \; ) \; + \varphi_1 z^n \; ( \; t- \frac{\tau}{2}
\; ) \; \end{equation}
where $z = ( {\bf u}, {\bf i} )$  represent finite integrals  over
electric and magnetic field strengths,  $\psi_0$ essentially interchanges 
${\bf u}$  and  ${\bf i}$ and $\varphi_0$, $\psi_0$ are selfadjoint $\mathbb{R}$-linear operators, cf. references \cite{Hein11},\cite{Hein12} and the appendix of this paper.

By saying that a TLM system ($\mathcal{C}, \mathcal{R}_{\zeta}$) generates solutions of equations 
(\ref {4.2}) it is meant that the total fields $z^{n,p}$, cf.(\ref {2.3},\ref {3.4}), 
satisfy the equations for {\it all}  vectors $z_{in}^p$  incident at any cell
(and as step functions of time of course matching $\tau$). Obviously, this
defines a  {\it  local}, so far purely algebraic property of the TLM system,
i.e. one referring to $\mathcal{R}_{\zeta}$ and indepedent of $\mathcal{C}$ that still
disregards all questions of convergence or computational stability, for
instance. As a rule, the TLM system solving (\ref {4.2}) depends strongly on
the time step and additional contraction properties must be checked on selecting
a particular stable system (cf. sect.5).

The deflection lemma establishes general recurrence relations between
solutions of  model equations that differ by a perturbation. 

Let a TLM system ($\mathcal{C}, \mathcal{R}_{\zeta}$) generate solutions of equation (\ref {4.2}) with any $\mathcal{F}$ of the form
\begin{equation}\label{4.5}
\mathcal{F}  [ z_+^n ]  [z^p  ]
= \sum_{\mu \in \mathbb{N} } 
 \varphi_{\mu} z^n \; ( \; t+ \frac {\tau}{2} - \mu \tau \; ) \; +
 \psi_{\mu} z^p \left(t - \mu \tau \right) 
\end{equation}
wherein $\varphi_{\mu}$, $\psi _{\mu}$ are $\mathbb{R}$-linear operators (many if
not almost all of which may be zero). Consider then a perturbation of
equation (\ref {4.2}) by a causal possibly non-linear operator $\mathcal{J}$ that
maps into the image  space of $\mathcal{F}$ (i.e.  $\varphi_{\mu}$, $\psi _{\mu}$
and $\mathcal{J}$ share the same image space). $\mathcal{F}$ and $\mathcal{J}$  may be time
dependent, i.e. step functions of time matching $\tau$.

Then one can state the  following 

\begin{uproposition}[Deflection Lemma]
A TLM system ($\mathcal{C}, \mathcal{R}_{\zeta}^{\sim}$),  with the same  connection map as  above, generates solutions of equation 
\begin{equation}\label{4.6} \mathcal{F}  [z_+^n ] [ z^p ] = \mathcal{J} [ z_-^n] [z^p ]
\end{equation} if and only if the {\it deflection}
\begin{equation}\label{4.7}
\mathcal{D} :=\mathcal{R}^{\sim} -\mathcal{R}
\end{equation}
satisfies recursively
\begin{equation}\label{4.8}
\varphi _{o,t} \mathcal{D}_t =
-\mathcal{J}_t [z_-^n ] \left[z^p \right]
 - \sum_{\mu \in \mathbb{N}_{+}}  \left( \varphi_{\mu} + \psi _{\mu -1} \right)_t
   \mathcal{D}_{t- \mu \tau} 
\end{equation}
\end{uproposition}

\begin{ucorollary}
Let $\varphi_0$ be injective (one to one) on the outgoing subspace and assume 
$\varphi _0 \circ \pi_{out}$ covers the ranges of $\mathcal{J}$, $ \varphi_\mu \circ
\pi_{out}$, and $\psi _{\mu} \circ \pi_{out}$; $\mu \in \mathbb{N}$. Then
\begin{equation}\label{4.9} 
\mathcal{D}_t := - \left(\varphi_0 \circ \pi_{out} \right)_t^{-1}
\; \{ \; \mathcal{J}_t [z_-^n] [z^p] + \sum_{\mu \in \mathbb{N}_+}
\left(\varphi_\mu + \psi_{\mu -1} \right)_t \mathcal{D}_{t- \mu \tau} \; \} \;
\end{equation}
defines recursively a causal operator $\mathcal{D}$ such that
\begin{center}
$\mathcal{R}^{\sim} := \mathcal{R} + \mathcal{D}$
\end{center}
is the reflection map of a TLM system ($\mathcal{C}$, $\mathcal{R}_\zeta ^{\sim}$)  which generates
solutions of equation (\ref{4.6}).

Moreover, if $\mathcal{D}_0 \equiv 0$, then this system is the {\it unique} one which
branches from  ($\mathcal{C}, \mathcal{R}_\zeta$) at $t=0$, i.e. that initially coincedes
with the unperturbed system. 
\end{ucorollary}

Condensed statements of that kind support a comment before being proved:

Guided by conceivable applications the perturbation $\mathcal{J}$ is written as a
function  of total port and node states and of time. Most commenly, it is a
function of nodal states only.

If $\varphi_0 \circ \pi_{out}$ is invertible,  then the negative time shift in
(\ref {4.6}) makes   $\mathcal{J}_t$ depend on nodal states up to
\begin{equation}\label{4.10}
 z^n \; ( \; t- \frac {\tau}{2} \; ) \; = z_{in}^p
\left( t- \tau \right) + \mathcal{R}^{\sim}  \left[ z_{in}^p \right]_{t- \tau}  =
z_{in}^p \left( t- \tau \right) + \mathcal{R}  \left[ z_{in}^p \right]_{t- \tau} + \mathcal{D}_{t-
\tau} \end{equation} 
thus ensuring that $\mathcal{D}$ enters the right-hand side of
recursion (\ref {4.9}) (implicitely in the argument of $\mathcal{J}_t$) up to $\mathcal{D}_{t-
\tau}$ at most, rather than up to $\mathcal{D}_t$ which obviously would lead to
ill-defined recurrence relations  without
further severe restrictions on $\mathcal{J}$. 

Note that from $\mathbb{R}$-linearity only additivity of $\varphi_{\mu}$, $\psi_{\mu}$ enters the following 
\begin{proof}
Let hence ($\mathcal{C}, \mathcal{R}_{\zeta}$) generate solutions of (\ref {4.2}) and let for a
second TLM system ($\mathcal{C}, \mathcal{R}_{\zeta}^{\sim}$)  the deflection $\mathcal{D}$ be defined
by (\ref {4.7}). Then for every sequence of states $[ z_{in}^p ]$
incident at a cell the total states of the second system are
\begin{eqnarray}\label{4.11}
&\tilde{z}^p\left(t \right) & =z_{in}^p +\mathcal{R}^{\sim} 
\left[z_{in}^p \right]_{t- \tau} \nonumber \\
& &=z_{in}^p +\mathcal{R} \left[z_{in}^p \right]_{t- \tau} 
+ \mathcal{D}_{t- \tau} \nonumber \\
& \tilde{z}^n \left( t + \frac {\tau}{2} \right) & = z_{in}^p 
+ \mathcal{R} \left[z_{in}^p \right]_{t} +\mathcal{D}_t 
\end{eqnarray}
This substituted for $z^{n,p}$ in equation (\ref {4.6}) yields
\begin{equation}\label{4.12}
\begin{array}{rcl}
\mathcal{F}[ \tilde{z}_+^p ]  [ \tilde{z}^p ] 
& = & 
\sum_{\mu \in \mathbb{N}}  \; \{ \; \varphi_\mu  ( z_{in}^p ( t- \mu \tau )   
+  \mathcal{R} [ z_{in}^p ] _{t- \mu \tau} 
+ \mathcal{D}_{t-  \mu \tau} ) \\
& \quad &   
 + \psi_\mu ( z_{in}^p ( t- \mu \tau )  +
\mathcal{R} [ z_{in}^p]_{ t-\tau- \mu \tau}  
+ \mathcal{D}_{ t- \tau - \mu \tau} )  \}   \\ 
& = & \sum_{\mu \in \mathbb{N}}  \; \{ \;  \varphi_\mu   
( z^n ( t + \frac {\tau}{2} -\mu\tau )  
+ \mathcal{D}_{t- \mu \tau} )   \\
& \quad & 
+ \psi_\mu ( z^p ( t -\mu \tau )  
+ \mathcal{D}_{t- \tau -\mu \tau} ) \}   \\     
& = & \mathcal{J}_t [ \tilde{z}_-^n ]  [ \tilde{z}^p ] 
\end{array}
\end{equation}
By virtue of the additivity of $\varphi_{\mu}$, $\chi_{\mu}$ follows
\begin{equation}\label{4.13}
\mathcal{F}  [ \tilde{z}_+^n ][ \tilde{z}^p ] = \mathcal{F}[z_+^n ] [ z^p ]   +
\sum_{\mu \in \mathbb{N}} \varphi_\mu \mathcal{D}_{t-\mu\tau} + \psi_\mu \mathcal{D}_{t-\tau -\mu\tau}=
\mathcal{J}_t[ \tilde{z}_-^n ] [ \tilde{z}^p]
\end{equation}

Equation (\ref {4.2}) states that $\mathcal{F}[z_+^n ] [ z^p ] $ above vanishes and inspection of the remaining terms yields the proposition.
\end{proof}

\vspace{-8pt}
\section {The S-matrix propagator}
The first order linear equation (\ref{4.4}) admits  a solution in closed form
as an iterated scattering process. This has been pointed out in detail in
references \cite{Hein11}, \cite{Hein12} at the example of the non-orthogonal
Maxwell field model. Using the commuting family of projections which split the
state space of each cell $\zeta$, with identity operator $Id_\zeta$, into
observable link states and unobservable gauge or stub vectors
\begin{equation}\label{5.1} Id_\zeta = \pi_l \oplus \pi_s
\end{equation}
and then again the link states into incoming and outgoing
subspaces of port and node  quantities 
\begin{equation}\label{5.2} 
\pi_l = \pi_{in} \oplus \pi_{out} 
\end{equation}
the reflection map of the first order linear process solving (\ref{4.4}) is
\begin{eqnarray}\label{5.3}  
z_{out}^p (t+ \tau)  &=&\mathcal{R} \left[z_{in}^p \right]\,
=\,z_{out}^n \; ( \; t+ \frac {\tau}{2} \; ) \;   \nonumber \\               
&=& \pi_{out}\; S \;z_{in}^n \; ( \; t+  \frac {\tau}{2} \; ) \;              
\nonumber \\  & &+ \; \pi_{out} \;
S \; \sum_{\mu
=1}^{\infty}  \left( \pi_s S \; \right)^\mu z_{in}^n 
\; ( \; t+ \frac {\tau}{2} - \mu \tau \; ) \;    
\end{eqnarray}
S denotes the nodal scattering operator which acts between in and outgoing
node  and stub vectors and can thus be written
\begin{equation}\label{5.4} 
S= \left(\pi_{out} \oplus \pi_s \right) S\left(
\pi_{in} \oplus \pi_s \right)=K+L+M+N  
\end{equation}
where as usual 
\begin{equation}\label{5.5} 
K:=\pi_{out}S\pi_{in},\quad
L:=\pi_{out}S\pi_{s},\quad M:=\pi_{s}S\pi_{in},\quad N:=\pi_{s}S\pi_{s}
\end{equation}  
have been introduced. With these operators equation (\ref{5.3}) 
reads more concisely 
\begin{equation}\label{5.6}
z_{out}^n \left( t \right) = K z_{in}^n \left( t \right) + L \sum _{\mu =
0}^{\infty} N^\mu M z_{in}^n \left( t- \left(\mu +1 \right) \tau \right)
\end{equation} 
(graphically represented in fig.2a).

$K$, $L$, $M$, $N$ are easily determined from equation (\ref{4.4})
by substituting for $ z^{n,p} = z_{in}^{n,p} +  z_{out}^{n,p}$ in
this equation the total scattering response of a Dirac pulse
\begin{equation}\label{5.7} z_{in}^p \left(t \right) = z_{in}^n ( \; t+
\frac{\tau}{2} \; ) \; =  \left\{  \begin{array}{rl}
z_0 &\;  if \qquad t \in [o, \tau) \\
0  &\;  else \quad 
\end{array}
\right. 
\end{equation}
applied and added to (\ref{5.3}, \ref{5.6}). By evaluating
equations (\ref{4.4})
 with these replacements at time intervals $[ k \tau , (k+1) \tau )$, $k \in
\mathbb{N}$, provides   operator identities that are independent for
$k=0,1$ and $k \ge 2$ under rather general conditions, cf. \cite{Hein9}
- \cite{Hein12} and which allow to determine $K$ uniquely, while $L$,$M$, 
and  $N$  are unique at most up to a linear automorphism of the stub space.
Indeed, sheer inspection of the scattering response (\ref {5.3}, \ref{5.6})
makes clear that any invertible linear transformation $G$ of the
stub space yields an equivalent scattering representation of $\mathcal{R}$
if $S$ is replaced by
\begin{equation}\label{5.8} 
S' := K+L'+M'+N' 
\end{equation} 
with $L':=LG^{-1}$, $M':=GM$,  $N':=GNG^{-1}$. It is this property that led 
to the understanding of  stub vectors as  internal {\it gauge} degrees of
freedom,  the automorphisms  of the stub linear space (in transmission
line parlance: a group of generalized impedance transformations)  constituting
the pertinent {\it gauge group}.

As already mentioned, the inherent gauge features of the TLM scheme can actually
be traced back to the deflection lemma, the  latter displaying in like
manner the generation mechanism for stub degrees of freedom as rather
suggestively also their gauge property.

Given a TLM system ($\mathcal{C}$, $\mathcal{R}_\zeta$) that generates 
solutions of equations (\ref{4.2}) with any operator $\varphi_0$
sharing the technical prerequisites refered to,
it is a straightforward exercise applying
the lemma and its corollary to write down the updating instructions for a
`deflected' system ($\mathcal{C},\mathcal{R}^\sim$) that solves the perturbed
equations (\ref{4.6}):
\begin{eqnarray}\label{5.9}
z_{out}^p \left(t + \tau \right) &:=& \mathcal{R}_t 
\left[ z_{in}^p 
\right]+ \mathcal{D}\left( t\right)\nonumber\\ 
I \left(t + \tau \right) &:=& - \left( \varphi_0 \circ \pi_{out} 
\right) _{t+ \tau}^{-1} \mathcal{J}_{t+ \tau}  [z_+^n][z^p] \\ 
\mathcal{D} \left(t + \tau \right) &:=& I \left( t + \tau \right)
- \left( \varphi_0 \circ \pi_{out} \right) _{t+ \tau}^{-1} 
\sum_{\mu \in \mathbb{N}} \left(\psi_\mu + \varphi_{\mu +1}
\right)_{t+ \tau}
\mathcal{D}\left( t- \mu \tau \right)\nonumber
\end{eqnarray} 
with all updating operations done in this order and implicitly 
after the first line 
\vspace{-5pt}
\begin{center}
$z^n \left(t+ \frac{\tau}{2} \right)
= z_{in}^p\left(t \right)+ z_{out}^p \left(t+ \tau \right)$ 
\end{center} 
Besides the proper interaction term $I$, which essentially appears in the model equations (\ref{4.6}) and is thus amenable to a direct physical interpretation, the deflecting field $\mathcal{D}$
introduces additional degrees of freedom that in general remain without
counterparts in the equations.

The gauge character of the deflection field is suggestively illustrated within
scattering representations (\ref{5.3},\ref{5.6}) for $\mathcal{R}$ and $\mathcal{R}^\sim$ -
assuming for once that such exist. If $S^\sim$ denotes the S-matrix of 
$\mathcal{R}^\sim$, then it admits a block decomposition in matrix form
\begin{equation}\label{5.10} S^\sim = \left( \begin{array}{ccc}          S    
      &\vline&       L^\sim   \\ \hline          M^\sim &\vline&       N^\sim 
\end{array}
 \right)
\end{equation}
wherein S represents the S-matrix of $\mathcal{R}$ and $L^\sim$, $N^\sim$ act on the additional 
deflection field $\mathcal{D}$.  The updating instructions (\ref{5.9}) 
obviously fit with
$L^\sim$ of the form 
\vspace{-5pt}
\begin{equation}\label{5.11} L^\sim = \left(
\begin{array}{ccc}
1 & \quad & \quad \\ 
\quad & \ddots & \quad \\
\quad & \quad& 1 \quad 0 \ldots 0
\end{array}
\right)^{T}
\end{equation}
(all other matrix elements of $L^\sim$ being zero). 
Any invertible transformation $G$, according to (\ref{5.8}) simultaneously applied to $L^\sim$, $M^\sim$ and $N^\sim$, then cannot alter $\mathcal{R}^\sim$ and thus
yields an equivalent representation. (The reader may trace this forth to  a
set of modified updating equations (\ref{5.9}).)

Spatial deformations of the TLM mesh continously get through stubs. In fact
they are modelled using stubs \cite{Chris1}, \cite{Hoef}, \cite{Hein9}. Therefore smooth
connections between mesh point dislocations and infinitesimal impedance
transformations $G$ exist that textually complete the gauge field analogy
in the sense of  fibrations \cite{May3}.

Once the gauge fixture of TLM being made available, full advantage can be taken of its
power and flexibility.  Aside from the nonlinear potentiality opened by the
deflection lemma, manifold variants of the linear scheme (\ref{5.6}) are
conceivable. For instance the generalized ansatz 
\vspace{-5pt}
\begin{equation}\label{5.12}
z_{out}^n (t)  = 
Kz_{in}^n (t) 
+ \sum_{\kappa =1}^k L_\kappa \sum_{\mu \in \mathbb{N}}
N_\kappa^\mu M_\kappa  z^n_{in} \left( t-\left(\mu + \kappa \right) \tau \right)
\end{equation} 
(which is graphically depicted in fig.2b) is suited to solve linear
model equations of the class (\ref{4.5}) with time differences in the port
states up to order~$k$.
\nopagebreak
\vspace{2.0cm}
\begin{figure}[ht]
\psset{xunit=1.0cm,yunit=1.0cm}
\begin{minipage}[t]{4.5cm} 
\begin{picture}(4.5,1.5)
\rput(2.35,0.95){\includegraphics[scale=0.7,
clip=0,bb=214 382 380 459]{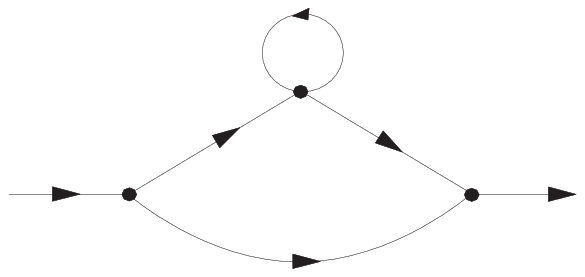}}
\rput(2.35,2.0){\makebox(0,0)[lb]{\scriptsize N}}
\rput(2.35,1.6){\makebox(0,0)[lb] {\small $\tau$}}
\rput(1.5,1.1){\makebox(0,0)[lb]{\scriptsize M}}
\rput(1.9,0.7){\makebox(0,0)[lb]{\small $\tau$}}
\rput(3.1,1.1){\makebox(0,0)[lb]{\scriptsize L}}
\rput(2.7,0.7){\makebox(0,0)[lb]{\small $\tau$}}
\rput(0.01,0.58){\makebox(0,0)[lc]{\small $z_{in}^p$}}
\rput(4.5,0.58){\makebox(0,0)[lc]{\small $z_{out}^p$}}
\rput(2.43,-0.23){\makebox(0,0)[lb]{\scriptsize K}}
\rput(2.33,0.3){\makebox(0,0)[lb]{\small $\tau$}}
\end{picture}
\begin{flushright}
a)
\end{flushright}
\end{minipage}
\hfill
\begin{minipage}[t]{6.0cm} 
\begin{picture}(5.8,1.5)
\rput(2.8,1.07){\includegraphics[scale=0.7,
clip=0,bb=192 377 402 463]{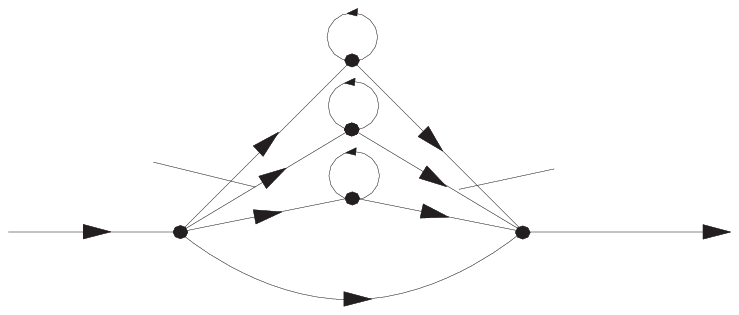}}
\rput(2.55,2.2){\makebox(0,0)[lb]{\scriptsize $N_\nu$}}
\rput(2.60,2.5){\makebox(0,0)[lb] {\tiny $\vdots $}}
\rput(1.7,1.25){\makebox(0,0)[lb]{\scriptsize $2\tau$}}
\rput(3.35,1.25){\makebox(0,0)[lb]{\scriptsize $2\tau$}}
\rput(0.78,0.95){\makebox(0,0)[lb]{\scriptsize $\frac{3\tau}{2}$}}
\rput(4.25,0.9){\makebox(0,0)[lb]{\scriptsize $\frac{3\tau}{2}$}}
\rput(2.05,0.45){\makebox(0,0)[lb]{\small $\tau$}}
\rput(3.15,0.45){\makebox(0,0)[lb]{\small $\tau$}}
\rput(-0.15,0.56){\makebox(0,0)[lc]{\small $z_{in}^p$}}
\rput(5.45,0.56){\makebox(0,0)[lc]{\small $z_{out}^p$}}
\rput(2.60,-0.27){\makebox(0,0)[lb]{\scriptsize K}}
\rput(2.60,0.25){\makebox(0,0)[lb]{\small $\tau$}}
\end{picture}
\begin{flushright}
b)
\end{flushright}
\end{minipage}
\begin{flushleft}
Figure 2: Linear scattering schemes \\
Every path from left to right following the arrows yields a delayed or phase
shifted contribution to the scattering response.
(Loops may be repeatedly run through any number of times.)
\end{flushleft}
\end{figure}

\vspace{-15pt}
A direct way to prove the above statement and to derive S-matrix
blocks $K$, $L_\kappa$, $M_\kappa$ and $N_\kappa$ in (\ref{5.12}) is again by
evaluating equation (\ref{4.5}) at successive time steps for an incident Dirac
pulse (\ref{5.7}).

\vspace{-8pt}
\section {The stable process in the time and frequency domain}
Parallel to the TLM typical splitting of states into link and gauge degrees of
freedom external and internal conditions for computational stability may be
distinguished, which respectively stress some physical aspects and a purely
numerical origin of instability.

In any interpretation the total link vectors are directly related to 
physical quantities that are in general subjugate to conservations laws. 
It is hence merely a strengthened consistence postulate to require that 
the reflection map of each cell must comply with these laws. 
The entire reflection response of an incident Dirac pulse integrated 
over time has for instance to be physically balanced with the input.

Actually this goes along with the Lax-Phillips contraction conditions 
\cite{Lax2} which are sufficient to prevent the observable fields from 
increasing indefinitely in the mesh. To implement these conditions numerically
in the TLM specifical fashion, internal algorithm stability has in addition
to prevent gauge fields (i.e. the 'stub quantities') from piling up 
to infinity. 
For instance, in a stable non-linear perturbed process which solves equation
(\ref{4.6}) the deflection must die out in time for any transitory incident
fields that do not pass over a certain level.

A universal condition that guarantees computational stability, even in the 
non-linear case (at least up to some maximum level of incidence), 
is quite evidently all updating operations must be contractive in a 
neighbourhood of zero. For the linear processes (\ref{5.6}), (\ref{5.12}) 
the operators $N_{(\kappa)}$ control updating between gauge fields. 
Hence a sufficient condition for internal stability of these processes 
is
\begin{equation}\label{6.1}
\parallel N \parallel < 1 
\end{equation}
if the norm $\parallel \ldots \parallel $ is submutiplicative, 
i.e.  characterized by
\begin{equation}\label{6.2}
\parallel N^\mu \parallel \; \le\; \parallel N \parallel ^\mu \; , 
\qquad \mu \in \mathbb{N}
\end{equation} 
The sup-norm with respect to any state space norm $\parallel z \parallel$
\begin{equation}\label{6.3}
\parallel N \parallel_{\sup} \; := \; \sup_{\parallel z \parallel =1}
\parallel Nz \parallel  
\vspace{-5pt}
\end{equation}
or the Hilbert (spectral) norm
\begin{eqnarray*}
\parallel  N \parallel_{H} \; &:=& \; \max \mid \lambda \mid \nonumber \\
& & \;{\lambda \;  eigenvalue \;  of \: N}
\end{eqnarray*} 
are, for instance, submultiplicative in the sense of (\ref{6.2}). 
All operator norms occurring in the following are tacitely assumed to share 
this property.

Condition (\ref{6.1}) is in general realized in the cases of interest 
by introducing bounds for the time step.
In the frequency domain, for a process of angular frequency 
$\omega =2 \pi f$, i.e.
\begin{equation}\label{6.5}
z_{in, out}^p \left(t+ \mu\tau \right)= 
e^{j\omega\mu\tau} z_{in,out}^p \left(t \right)\;,\qquad \mu \in \mathbb{Z}
\end{equation}
equation (\ref{5.6}) reads
\begin{equation}\label{6.6}
z_{out}^n = \left( K+L \sum_{\mu=0}^\infty N^\mu M e^{j \omega \left( \mu + 1 \right) \tau} \right) z_{in}^n
\end{equation}
With $\mathbb{C}$-linear operators $K$, $L$, $M$, $N$ this turns into
\begin{equation}\label{6.7}
z_{out}^n = \left(K+e^{-j\omega\tau} L \sum_{\mu=0}^{\infty}
\left( e^{-j\omega\tau}N  \right)^\mu M  \right) z_{in}^n
\end{equation}
A stable process, characterized by condition (\ref{6.1}), ensures that the power series in $e^{-j\omega\tau}N$ converges to $(Id-e^{-j\omega\tau}N)^{-1}$. Equations (\ref{6.7})  thus becomes $z_{out}^n = S^\sim z_{in}^n$ with
\begin{equation}\label{6.8}
S^\sim = K+L\left(e^{j\omega\tau} Id-N  \right)^{-1} M
\end{equation}
Evidently, in the frequency domain the complex S-matrix connects directly incident with outgoing fields and stubs do no longer enter the algorithm in the form of explicitely computed quantities. Instead, the stubs yield implicitely an additive contribution to the S-matrix 
\begin{equation}\label{6.9}
S_g^\sim =L\left( e^{j\omega\tau} Id-N  \right)^{-1} M
\end{equation}
extracted from decomposition (\ref{6.8}). It is virtually this gauge contribution that anew characterizes the TLM scheme here in the frequency domain.

Of course, knowing the nodal S-parameters of each cell, 
'intrinsic S-matrices' of entire mesh systems can be computed, 
just as then S-parameters for any combination of submeshes. 
(Proceeding in that direction opens 
rather a special field of application to linear network theory 
than deeper structural insight into the TLM method).

In the frequency domain the TLM algorithm solves the eigenvector equation
\vspace{-7pt}
\begin{equation}\label{6.10}
\left(\mathcal{C} S_\zeta^\sim e^{j\omega\tau} - Id \right)\left( z_{in} 
\right)_\zeta = 0
\end{equation}
where $\zeta$ runs over all mesh cells and the connection map $\mathcal{C}$ interlaces a steady state excitation which may be imposed on the mesh boundary according to any mode templates for instance.

Various methods for solving equations (\ref{6.10}) are imaginable.
We experimented with classical conjugate gradient schemes, such as 
Fletcher-Reeves and Polak-\-Ribi\`ere which, however, showed poor 
convergence for larger mesh systems of several thousand cells. 
Very satisfactory results have been obtained, in contrast, with  
a fixed point iteration analogous to a static time domain process
(of course with complex fields). 
Indeed, reiterating the instructions
\vspace{-7pt}
\begin{eqnarray}\centering\label{6.11}
\left( z_{in} \right)_{\zeta , k} & = & \mathcal{C}
\left( z_{out} + z_{exc} \right)_{\zeta, k} \nonumber \\
\left( z_{out} \right)_{\zeta , k+1} & = & e^{j\omega\tau}
S_\zeta^\sim \left( z_{in} \right)_{\zeta, k}
\end{eqnarray}
converges for any steady state excitation $(z_{exc})_\zeta$ quite rapidly 
to a fixed point, and thus to a solution of equation (\ref{6.10}) which 
matches the imposed excitation pattern. Still faster convergence had been
attained by starting with zero fields and smoothly approaching the steady
state excitation over a transitory phase of some hundred iterations.

Fig.3 displays convergence thus attained for a 3 dB waveguide 
directional Riblet coupler.
\vspace{-8pt}
\input{f3.tex}
\vspace{-8pt}
\section {Charged particles coupled to the Maxwell field}
It has been stressed that this study focusses on the intrinsic 
structure of the TLM algorithm, not on any interpretation. 
Still one application is sketched in this section to illustrate 
the excellent capacity of TLM to cooperate with a tied finite 
difference scheme. Beyond that, a variety of original research 
applies concepts of this paper and should be consulted for technical
illustration. In particular, part of references \cite{Hein8} - \cite{Hein12}
constitute by now `standard technicality' which can be referred to in 
the following. (A bridge between the non-orthogonal Maxwell field model, 
as treated in \cite{Hein12}, and the presently generalized formalism is 
built in the appendix.)

The following application extends the Maxwell field model in a 
non-linear fashion. By deflection, a relativistic charged particle 
current is coupled to the apriori free electromagnetic field. 
The underlying dynamical relations are thus Maxwell's equations 
and the equations of motion for a particle current with interactions 
being brought in by Amp\`ere's law and the Lorentz force.

The relativistic mass of a particle at velocity ${\bf v}$ is
\begin{equation}\label{7.1}
m=m_0 \; ( \; 1- \frac{v^2}{c_0^2} \; )^{\frac{-1}{2}}
\end{equation}
where $m_0$ and $c_0$ respectively denote the particle rest mass and
the velocity of light.  For particles of charge $q_0$ the Lorentz force 
in the laboratory frame is 
\vspace{-5pt}
\begin{equation}\label{7.2}
{\bf F}_L = \frac{d}{dt} \left( m {\bf v}  \right) =q_0 \left( {\bf E}+{\bf v }\wedge {\bf B} \right) -\nu_c m {\bf v} 
\end{equation}
where a mean collision frequency $\nu_c$ has been introduced to take internal friction into account; cf. \cite{Cap}.\newline
Equation (\ref{7.2}) is linear in $d{\bf v}/dt$ and reads explicitely
\begin{equation}\label{7.3}
\frac{d}{dt}\left(m{\bf v} \right) = m \frac{d{\bf v}}{dt}+{\bf v}\frac{dm}{dt}=\mathcal {M}\frac{d{\bf v}}{dt}={\bf F}_L
\end{equation}
with
\vspace{-5pt}
\begin{equation}\label{7.4}
\mathcal {M}=m \left(\delta_{ij} + \lambda v_i v_j  \right)_{i,j=1,2,3} \quad ; \quad \lambda =\left(c_0^2 -v^2 \right)^{-1} \end{equation}
by virtue of (\ref{7.1}). $\mathcal {M}$ is obviously selfadjoint and 
strictly positive for $v < c_0$, hence invertible, and yields the 
relativistic particle acceleration as 
\vspace{-4pt}
\begin{equation}\label{7.5}
\dot {\bf  v} =\mathcal {M}^{-1} {\bf F}_L 
\end{equation}
Then $n$ particles per unit volume of (mean) velocity ${\bf v}$ generate a current density
\begin{equation}\label{7.6}
{\bf j}_c = \rho \cdot {\bf v} \quad ; \quad \rho = nq_0
\end{equation}
which has to be included in the generalized Amp\`ere's law (1st Maxwell's integral equation)
\begin{equation}\label{7.7}
\int\limits_{\partial A} {\bf H}d{\bf s} =
\int\limits_{A} {\bf \varepsilon} \frac{\partial {\bf E}}{\partial t} d{\bf S} +
\int\limits_{A} {\bf \kappa}_e {\bf E} d{\bf S} +
\int\limits_{A} {\bf j}_c d{\bf S}
\end{equation}
Charge conservation is ensured by the continuity equation for $( \rho, {\bf j}_c)$ or its integral equivalent, the Gauss-Ostrogadski law
\begin{equation}\label{7.8}
\int\limits_{\partial C}{\bf j}_c d {\bf S} = - \frac{d}{dt} \int\limits_C \rho d V 
\end{equation}
valid for arbitrary test volumes $C$. Of course Faraday's law (Maxwell's 2nd integral equation) holds unalteredly 
\begin{equation}\label{7.9}
- \int\limits_{\partial A} {\bf E}d{\bf s} =
\int\limits_{A} {\bf \mu} \frac{\partial {\bf H}}{\partial t} d{\bf S} +
\int\limits_{A} {\bf \kappa}_m {\bf H} d{\bf S}
\end{equation}
Perhaps, in the present context the loss currents called forth by the electric and magnetic conductivities ${\bf \kappa}_e$, ${\bf \kappa}_m$ in Maxwell's equations (\ref{7.7}), (\ref{7.9})  attract some attention. They simply allow for simulating dissipative effects of any origin in the plasma. Beyond it, further simplifying assumptions obviously enter the above description, which may be modified or dropped in each case. Thus, the model neglects the direct repulsive Coulomb interaction between the particles, which amounts to a low density approximation. Also, the fixed mean collision frequency $\nu_c$  may be substituted by a state dependent expression.

In the TLM model the physical state, given by the Maxwell field $({\bf E}, {\bf H})$, the charge density $\rho$, and the mean particle velocity ${\bf v}$ has to be related to cell states in a suitable interpretation. In a convex hexahedral parcel twines cell (wherein ports on the cell faces interconnect the midpoints of the cell edges, cf. fig.1) the fields are connected to cell states in terms of port voltages and currents
\begin{equation}\label{7.10}
U_\mu := {\bf E} \cdot {\bf p}^\mu \quad , \quad I_\mu := +(-) {\bf H} \cdot {\bf p}^\nu
\end{equation}
the sign corresponding with that in \begin{equation}\label{7.11}
{\bf f} = +(-) {\bf p}^\mu \wedge {\bf p}^\nu
\end{equation}
where ${\bf f}$ denotes the cell face vector which points into the cell.
On the cell surface $U_\mu^p$, $I_\mu^p$ as finite integrals approximate line integrals
\begin{equation}\label{7.12}
U_\mu^p \approx \int\limits_{p_\mu}{\bf E} d{\bf s} \quad ; \quad 
I_\mu^p \approx +(-) \int\limits_{p_\nu}{\bf H} d{\bf s} 
\end{equation}
(every port vector in parcel-twines location is naturally identified with an integration path on the pertinent face).  

$U_\mu^n$, $I_\mu^n$  measure  the field strengths {\it in} the cell according to (\ref{7.10}),   yet with the nodal images ${\bf p}^{\mu,\nu}$, cf. sect.2, (\ref{2.1}, \ref{2.2}). The model equations for Amp\`ere's and Faraday's laws are developed in detail in references \cite{Hein9}, \cite{Hein11},  \cite{Hein12}.  A suitable change of representation, cf.\cite{Hein12} sect.4, provides decoupled equations in terms of transformed voltages and currents $( u_\mu^{n,p},i_\mu^{n,p})_{\mu=1,...,12}$. On the first components $z =(u_\mu , i_\mu )_{\mu =1,...,3}$ the discretized 1st Maxwell equation reduces to 
\begin{equation}\label{7.13}
\psi_0 z^p \left( t \right) = \varphi_0 z^n ( \; t + \frac{\tau}{2} \; ) \;
+ \varphi_1 z^n ( \; t - \frac{\tau}{2} \; ) \; +
{\bf J}_c
\end{equation}
the convection current ${\bf J}_c$ being added here as the following perturbation
\begin{equation}\label{7.14}
{\bf J}_c = \frac{1}{4}\det (B)B^{-1} \frac{Q}{V}{\bf v}= \frac{1}{4}B^{-1}Q{\bf v} 
\end{equation}
$Q$ is the total charge in the cell and ${\bf v}$ the mean particle 
velocity. The node vector matrix $B$ with adjoint $B_*$ is defined in 
reference \cite{Hein12} and (modulo a factor 4) also in reference 
\cite{Hein11}. In the present normalization the determinant of $B$ 
equals the cell volume and $A= \det (B)B_*^{-1}$ is the so called area 
matrix, cf.\cite{Hein12} sect.3. \newline
Still $Q=\rho V$ and ${\bf v}$ have to be incorporated into the TLM model 
as cell state functions. This is achieved by introducing six additional ports
in the form of the cell face vectors ${\bf f}^\mu$ which measure the 
convection currents through the respective cell boundary faces and through 
pertinent image areas in the node.\nopagebreak
The link states associated with these ports are 
\vspace{-3pt}
\begin{equation}\label{7.15}
{\bf J}_\mu := 
\left(
\begin{array}{rcl}
J_{in}& + & J_{out} \\
J_{in}& - & J_{out}
\end{array}
\right)
\end{equation}
to which the following interpretation is construed
\begin{equation}\label{7.16}
J_\mu := {\bf J}_{\mu,1} = J_{\mu,in} + J_{\mu, out} = \rho {\bf v} \cdot {\bf f}^\mu
\end{equation}
Note that voltages are not defined in parallel with these currents. 
This amounts to neglect a pressure in the here adopted low density i
approximation. In transmission line parlance, the characteristic impedance i
$\mathfrak{z}$ of the cell face channels is zero, just as then are  
the voltages $U_n= \mathfrak{z}J_{n,2}$, to which a pressure 
$p= \rho U$ is assigned. The generalized TLM scheme, which in the setting 
of sect.3 dispenses with impedances, still works for 
$\mathfrak{z}=0$  with link states defined as in (\ref{7.15}).
The total charge Q in a cell is related to the cell boundary current 
by charge conservation cf. (\ref{7.5}), which yields the updating 
instructions 
\vspace{-10pt}
\begin{equation}\label{7.17}
Q \left( t+\tau \right) = Q(t) + \tau \sum_{\mu=1}^6 J_\mu ^p (t)
\end{equation} 
In the cited canonical representation, \cite{Hein12} sect.4, 
the nodal states $z_k^n =(u_k^n , i_k^n)$ are interpreted as field strengths
\begin{equation}\label{7.18}
{\bf E} = B_*^{-1} \left( u_k^n \right)_{k=1,2,3} \quad;\quad {\bf H} = B_*^{-1} \left( i_k^n \right)_{k=7,8,9} 
\end{equation}
These substituted in equation (\ref{7.2}) yield the Lorentz force and updated particle velocities in the cells
\begin{equation}\label{7.19}
{\bf v}(t+\tau)={\bf v}(t)+\tau \mathcal {M}^{-1}{\bf F}_L
\end{equation}
for $\mathcal {M}({\bf v})$, ${\bf F}_L ({\bf v},{\bf E},{\bf H})$, and ${\bf v}$ at time $t$.

Ultimately, equations (\ref{7.13}), (\ref{7.16}), (\ref{7.17}), 
(\ref{7.19}) are simultaneously solved by deflection, branching out 
the trivial process $J^n(t)\equiv 0$ and the reflection map 
$\mathcal{R}$ of the free Maxwell field model as outlined in sect.4.
This results in the following updating instructions for the coupled
TLM process
\begin{eqnarray}
z_{out}^p(t+ \tau) & := & z_{in}^p(t) +\mathcal{R}[ z_{in} ^p]
+ \mathcal{D}(t) \\
q & :=  &  Q(t) + \tau \sum_{\mu=1}^6 J_{\mu ,in}^p(t) \\
{\bf v}(t+\tau) & := & {\bf v}(t) + \tau \mathcal {M}^{-1}{\bf F}_L \\
{\bf J}_c & := & \frac{1}{4}B^{-1} q {\bf v} \\
\mathcal{D}(t+\tau ) & := & \left( \pi_{out} \circ \varphi_0 
\right)^{-1} \left\{ \left(\psi_0 - \varphi_1 \right) \mathcal{D}(t)
+ {\bf J}_c \right\} \\
J_{\mu ,out}^p(t+\tau) & := & -J_{\mu ,in}^p (t) + \det (B)^{-1}q
{\bf v} \cdot {\bf f}^\mu \\
Q(t+ \tau) & := & q+ \tau \sum_{\mu =1}^6 J_{\mu ,out}^p(t + \tau)
\end{eqnarray}
with all operations carried out in this order.

Stability requires the assignments must be contractive in the domain of operation ( i.e. up to any allowed maximum excitation). This is ensured  with bounds for the time step following the guide-lines of sect.6.

\vspace{-8pt}
\begin {appendix}
\section {Appendix}
To remain technically in contact with the familiar framework 
(which uses line impedan\-ces) the bridge is thrown in this passage 
from Maxwell's equations in the transmission line setting \cite{Hein9},
\cite{Hein11}, \cite{Hein12} to the presently generalized formulation 
that works with projections instead.
In the canonical representation of the link state space, 
\cite{Hein12} sect.3, with 
\begin{equation}\label{app1}
z = {{\bf u} \choose {\bf i}} =\left(
\begin{array}{c}
{\bf u}_{in}+{\bf u}_{out} \\ 
y({\bf u}_{in} - {\bf u}_{out})      
\end{array} \right)
\end{equation} 
($y= 1 / \mathfrak{z}$ the characteristic field admittance) the 
generalized Amp\`ere's law discretized in a parcel-twines cell takes 
the form of equation (\ref{4.4}) with
\begin{equation}\label{app2}
\psi_0= \left(\begin{array}{cc}
0_3 & 1_3 \\ 
0_3 & 0_3 
\end{array} \right)
\end{equation}
($0_n, 1_n$ denote respectively the $n \times n$ zero and unit 
matrix blocks) and
\begin{equation}\label{app3}
\varphi_0 = \left(\begin{array}{cc}
T_+ & 0_3 \\
0_3 & 0_3 
\end{array} \right)\quad , \quad 
\varphi_1 = \left(\begin{array}{cc}
T_- & 0_3 \\ 
0_3 & 0_3 
\end{array} \right)\
\end{equation}
with real selfadjoint (i.e. symmetric) $3 \times 3$ matrix operators
\begin{equation}\label{app4}
T_{+(-)} = \frac{1}{4}\det (B)B^{-1} \left(\frac{ \kappa_e}{2} + (-) \frac{\varepsilon}{\tau} \right)B_*^{-1}
\end{equation}
B is the node vector matrix, defined ibid.,  which depends only on the cell geometry. $\kappa_e$ and $\varepsilon$ denote respectively the electric conductivity and permittivity tensors.
With (\ref{app1}) the projections into the incoming and outgoing states related to Amp\`ere's law (separated) are
\begin{equation}\label{app5}
\pi_{in}= \frac{1}{2}  \left(\begin{array}{cc}
1_3 & \mathfrak{z}1_3 \\ 
y1_3 & 1_3 
\end{array}\right)\quad , \quad
\pi_{out}= \frac{1}{2}\left( \begin{array}{cc}
1_3 & -\mathfrak{z}1_3 \\
-y1_3 & 1_3 
\end{array}\right)
\end{equation}
Proceeding as outlined in sect.5, i.e. substituting for $z^{n,p}$ in equation (\ref{4.4}) with operator coefficients (\ref{app2}), (\ref{app3}) the total scattering response of a Dirac pulse, cf. (\ref{5.6}), (\ref{5.7}), yields the S-matrix blocks $K_A,\ldots ,N_A$ pertinent to the generalized Amp\`ere's law. Modulo gauge transformations (\ref{5.6}) they are uniquely determined as\nopagebreak
\vspace{-5pt}
\begin{equation}\label{app6}
\begin{array}{lcl}
K_A & = & \pi_{out} U \pi_{in} \\
N_A & = & \pi_s V \pi_s \\
L_A & = & \pi_{out} W \left( Id - N_A^2\right)^{\frac{1}{2}} \pi_s \\
M_A & = & \pi_s \left( Id - N_A^2\right)^{\frac{1}{2}}W_* \pi_{in}
\end{array}
\end{equation}
with
\vspace{-5pt}
\begin{equation}\label{app7}
U= \left( \begin{array}{ccc}
yT_+^{-1} -1_3  & \vline & 0_3 \\
\hline    
0_3  & \vline & - yT_+^{-1} +1_3  
\end{array}
\right) \quad , \quad
V= - \; yT_+^{-1} - T_+^{-1} T_-
\end{equation}
W maps the stub vectors pertinent to Amp\`ere's law back into the link 
channels and acts in the canonical representation on $z_{l}^{n} \oplus 
z_{s}$ as a matrix operator
\begin{equation}\label{app8}
W=\pi^n_l W \pi_s = 2 \left( \begin{array}{ccc}  0
& \vline & 1_3 \\ \hline 0 &\vline& 0
\end{array} \right)
\end{equation} 
\noindent It is left as an exercise to show that $S_A^\sim$ with these data is unitary iff $\kappa_e=0$.

For Faraday's law the derivation of operators $K_F \ldots N_F$ 
runs completely parallel.
\end{appendix} 
\vspace{-8pt}
\section{Conclusions}

Despite the rich potentiality they offer, the gauge properties of the TLM
scheme remain still largely unexploited, today, in numerical model design.
A general method that allows for incorporating a linear or non-linear
interaction into a given linear model using these properties has been
developed in this paper,
on the algorithmic level and within a prototype application.

From a fundamental point of view, future work should yield still deeper insight
into the precise nature and the range of the observed gauge similarity.
It is clear from the preceeding that the latter is not to be (mis)understood in
the restricted sense of Lagrangian field theory: Until now, a {\it Lagrangian}
is not known in the TLM context (though it may be challenging to look for such
an object).
Instead, the structural relationship traced out
points to some more general common aspects: Much alike gauge field models,
the TLM algorithm works with unobservable quantities that intermediate
interactions and which together with their internal symmetries can be similarly
described in fibrations. There is a considerable differential
geometric overlap with gauge field theory that attracts our attention
and further study in this direction may significantly
promote the TLM method.

Ultimately, not a few technical spade-work of course remains to be done in any
particular application that is conceivable. The propagator approach
marks only a canonical path equipped with reliable instruments.


\hrulefill
\newline
\textsc{Steffen Hein};\;
DE-83043 Bad Aibling, Germany
\newline E-mail address:\; steffen.hein@bnro.de
\end{document}